\documentclass[draftclsnofoot, 12pt]{IEEEtran}

\usepackage{graphicx, epstopdf}
\usepackage{amsmath, amssymb, latexsym}
\usepackage{color}
\usepackage{multirow}
\usepackage{array}
\usepackage{hhline}

\begin{document}

\title{Analysis of the Packet Loss Probability in Energy Harvesting Cognitive Radio Networks}

\author{Shanai Wu,~\IEEEmembership{Student Member,~IEEE,}
        Yoan Shin$^\dagger$,~\IEEEmembership{Senior Member,~IEEE,} \\
        Jin Young Kim,~\IEEEmembership{Senior Member,~IEEE,}
        and Dong In Kim,~\IEEEmembership{Senior Member,~IEEE}
\thanks{This technical report is a detailed version of a paper which was accepted to appear in IEEE Communications Letters. Copyright belongs to IEEE.}
\thanks{This work was supported by the National Research Foundation of Korea (NRF) grant funded by the Korean government (MSIP) (2014R1A5A1011478).}
\thanks{Shanai Wu and Yoan Shin are with the School of Eletrical Engineering, Soongsil University, Seoul 06978, Korea (e-mail: \{sunae0814, yashin\}@ssu.ac.kr).}
\thanks{Jin Young Kim is with the Department of Wireless Communications, Kwangwoon University, Seoul 01897, Korea (e-mail: jinyoung@kw.ac.kr).}
\thanks{Dong In Kim is with the School of Information and Communication Engineering, Sungkyunkwan University, Suwon, Gyeonggi-do 16419, Korea (email: dikim@skku.ac.kr).}
\thanks{$^\dagger$ Corresponding author}
}
\markboth{[Technical Report]}%
{Shell \MakeLowercase{\textit{et al.}}: Bare Demo of IEEEtran.cls for Journals}

\maketitle

\begin{abstract}
A Markovian battery model is proposed to provide the variation of energy states for energy harvesting (EH) secondary users (SUs) in the EH cognitive radio network (CRN). Based on the proposed battery model, we derive the packet loss probability in the EH SUs due to sensing inaccuracy and energy outage. With the proposed analysis, the packet loss probability can easily be predicted and utilized to optimize the transmission policy (i.e., opportunities for successful transmission and EH) of EH SUs to improve their throughput. Especially, the proposed method can be applied to upper layer (scheduling and routing) optimization. To this end, we validate the proposed analysis through Monte-Carlo simulation and show an agreement between the analysis and simulation results.
\end{abstract}


\section{Introduction}
Energy harvesting (EH) allows devices to charge energy from ambient power sources, such as vibration, heat and electromagnetic waves. With the prevalence of wireless signals, radio frequency (RF) EH has received substantial attention and has dramatically grown. A comprehensive literature review on the research progress in wireless networks with RF EH is presented in [1]. RF EH becomes a candidate solution for self-sustainable energy supply, and the applications of RF EH are investigated in various systems, such as wireless charging systems [2] and wireless sensor networks [3]. It is particularly important for system such as body area networks (BANs), because the battery of such devices is implanted in human body and battery replacement is infeasible.

Besides energy efficiency, spectrum efficiency is one of the critical issues when designing a wireless network. To improve spectrum utilization efficiency, dynamic spectrum access (DSA) that allows the reuse of a radio spectrum in an opportunistic manner, has been proposed. The cognitive radio (CR) is the enabling technology to support DSA. In CR networks (CRNs), primary users (PUs) have higher priority in using radio spectrum. However, secondary users (SUs) are allowed to access the spectrum as long as the spectrum is not temporally used by the PUs, and share the spectrum as long as the PUs are properly protected [4]. Consequently, the EH CRN becomes a promising approach to support both spectrum efficiency and energy efficiency [5].

A considerable number of works focused on the EH CRNs from different perspectives. The opportunistic channel access policy is proposed to maximize the throughput of the EH SUs, when the SU accesses the spectrum to transmit data or to harvests energy if the spectrum is idle or occupied by the PU, respectively [6]. The optimization is formulated based on a Markov decision policy (MDP) and an online learning algorithm is applied to adapt the channel access policy without any prior knowledge. The energy causality constraint states that the total consumed energy should be equal to or less than the total harvested energy, and an optimal sensing policy is proposed to maximize the throughput of the EH SUs subject to an energy causality constraint and a collision constraint mandating that the PU is required to be properly protected [7]. In addition, to maximize the average throughput of the EH SUs, authors in [8] studied the optimal pair of the sensing duration and the sensing threshold under an energy causality constraint and a collision constraint. Since optimal sensing duration and sensing threshold are jointly intertwined with the sensing performance and energy causality constraint, it is difficult to directly acquire optimal sensing duration and sensing threshold. Therefore, the authors first obtained an optimal sensing duration by analyzing the effect of sensing duration for a given sensing threshold, and found an appropriate sensing threshold that meets the requirements of the previously derived optimal sensing duration. In [9], the impact of temporal correlation of the primary traffic on the achievable throughput of the EH SU is investigated. The temporal correlation of the primary traffic is modeled according to a time-homogeneous Markov process to derive the upper bounds on the achievable throughput as a function of the energy arrival rate and the temporal correlation of the primary traffic. Then the optimal sensing threshold is derived to maximize the theoretical upper bound on the achievable throughput under the energy causality and collision constraints. In [10], multiple SUs cooperate to sense a set of common channels, so as to maximize the overall EH CRN throughput. To address the cooperative optimization problem, an online learning algorithm is envolved to observe the primary channels and to make optimal decision when the SUs cooperate in a round-robin fashion. Moreover, the cooperative decentralized optimization is formulated as a decentralized partially observable MDP when the SUs cooperate in decentralized manner. According to [11], the advanced smart antenna technologies including multiple-input multiple-output and relaying techniques can be applied to significantly improve the energy efficiency and spectrum efficiency compared to single-antenna systems, and the corresponding receiver architecture is discussed.

However, a theoretical model that provide the variation of energy states in the EH SU was not developed. In this paper, we present a Markovian battery model for the EH SU equips with an RF transceiver and an RF energy harvester separately. From the proposed model, we derive the probability that the EH SU completely runs out of energy first. In the EH CRN, the EH SUs are allowed to transmit packets on a spectrum when the spectrum is idle and therefore may fail to transmit data in the absence of an idle spectrum or/and sufficient energy. Then, the packet loss probability in the EH SU due to sensing inaccuracy and energy outage is analyzed.

\section{Netwok Model}
An EH CRN with a single EH SU and numbers of PUs is considered, and a discrete time model with time slotted in intervals of length $T\in \mathcal{R}^{+}$ is assumed, where $\mathcal{R}^{+}$ denotes the set of nonnegative real numbers. In our work, the idle and the occupied spectrum states are denoted as $\mathcal {H}_{0}$ and $\mathcal {H}_{1}$, respectively. The spectrum state is changed with the generation of primary signals and modeled as a correlated, two-state process, as illustrated in Fig. 1. Therefore, the spectrum switches its states between idle and occupied randomly. If the spectrum is occupied by the PU in the current slot, the spectrum will be occupied by the PU with probability $q_{o}$ or the spectrum state will transit to idle with probability $1-q_{o}$ in the next time slot. Similarly, if the spectrum is idle in the current slot, then the spectrum will remain idle with probability $q_{i}$ or will be occupied by the PU with probability $1-q_{i}$ in the next slot. Then we can derive the steady state probability that the spectrum is idle or occupied as follows.
\begin{equation}
\left[\begin{array}{cc}
\pi_i & \pi_o \\ \end{array}\right] =
\left[\begin{array}{cc}
\pi_i & \pi_o \\ \end{array}\right]
\left[\begin{array}{cc}
q_i & 1-q_i \\
1-q_o & q_o \\ \end{array}\right]
\end{equation}

Consequently, the steady state probability of idle spectrum is derived as
\begin{equation}
\pi_{i}(q_{i}, q_{o})=\mbox{Pr}(\mathcal{H}_{0})=\frac{1-q_{o}} {2-q_{i}-q_{o}},
\end{equation}
and the steady state probability that the spectrum is occupied by the PU is $\pi_{o}(q_{i}, q_{o})=\mbox{Pr}(\mathcal{H}_{1})=1-\pi_{i}(q_{i}, q_{o})$.

\begin{figure}[!t]
\centering
\includegraphics[width=4in]{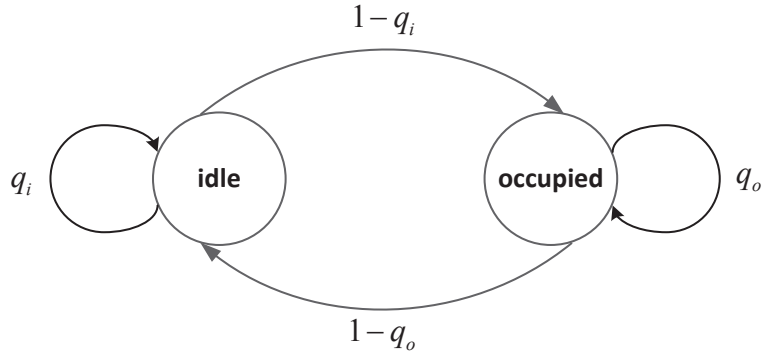}
\caption{Transition diagram for the spectrum states}
\label{Fig_sim}
\end{figure}

Suppose that the EH SU has an RF transceiver and an RF energy harvester (i.e., EH board) with separated antennas. Therefore, simultaneous data communication and EH are allowed in the EH SU. In other words, both {\it{in-band}} RF EH and {\it{out-of-band}} RF EH are supported with the general architecture introduced in [1]. While the EH SU can harvest energy from the same spectrum as that of data communication in the in-band RF EH, the process of EH board is assumed to be actively controllable so as not to share the same spectrum and reduce the efficiency of each process. Consequently, in this paper we assume out-of-band RF EH only to avoid the conflict between RF EH and data transmission, so that the EH SU with separated antenna harvests energy in the frequency spectrum different from that used for data communication. When the EH board is activated, the energy arrival is also modeled as a correlated, two-state process. Therefore, the EH board harvests energy at a constant rate $e_h$ per time slot as long as RF signals are detected; otherwise energy is not harvested. If the EH board harvests energy in the current slot, energy will be harvested in the next time slot with probability $p_{on}$ or will not be harvested in the next slot with probability $1-p_{on}$. However, if no energy is harvested in the current slot, no energy will be harvested with probability $p_{off}$ or energy will be harvested with probability $1-p_{off}$ in the next slot. Similarly, the steady state probabilities of energy harvested and no energy harvested are $e_{on}(p_{on}, p_{off})=\frac{1-p_{off}}{2-p_{on}-p_{off}}$ and $e_{off}(p_{on}, p_{off})=\frac{1-p_{on}}{2-p_{on}-p_{off}}$, respectively.

Under hypothesis $\mathcal{H}_{0}$ and $\mathcal{H}_{1}$, the discrete received signal at the EH SU can be given as 
\begin{equation}
y_{m}(n)=\left\{\begin{array}{lcl}
s(n)+w(n) & \mbox{for} & \mathcal{H}_{1}\\
w(n) & \mbox{for} & \mathcal{H}_{0}
\end{array}\right.
\end{equation}
where $y_{m}(n)$ represents $n$-th sample at the slot $m$. Both primary signal $s(n)$ and noise $w(n)$ are assumed to be modeled as circularly symmetric complex Gaussian (CSCG) with the variances $\sigma_{s}^2$ and $\sigma_{w}^2$, respectively [12].

The EH SU is assumed to have no priori information of primary channels. To opportunistically access the spectrum, the EH SU needs to perform spectrum sensing during sensing duration $\tau_{s}\in(0, T]$ at the beginning of each slot. The number of samples at the SU's energy detector is $N=\tau_{s}f_{s}$, when the received signal is sampled at a sampling frequency $f_{s}$. 
Then the test statistic for energy detection is given as
\begin{equation}
T(y_{m})=\frac{1}{N} \sum_{n=1}^N |y_{m}(n)|^2.
\end{equation}

We denote $\theta_{m}\in\{0(T(y_{m})<\varepsilon), 1(T(y_{m})>\varepsilon)\}$ as the spectrum sensing result at slot $m$. If the energy of detected signal at slot $m$ is larger than detection threshold $\varepsilon$, the spectrum is regarded as occupied by the PU and the sensing result $\theta_{m}$ is equal to 1; otherwise, the spectrum is idle and $\theta_{m}=0$. There are two important parameters associated with spectrum sensing: the false alarm probability and the detection probability.

Under the hypothesis $\mathcal{H}_{0}$, we denote $p_{0}(x)$ as the probability density function (PDF) of $T(y)$. For a chosen detection threshold $\varepsilon$, the false alarm probability $P_f$ that the presence of primary signal is falsely declared is given by
\begin{equation}
\begin{aligned}
P_{f}(\varepsilon) & \triangleq \mbox{Pr}\left(\theta_m=1|\mathcal{H}_{0}\right) \\
& = \mbox{Pr}\left(T(y_{m})>\varepsilon|\mathcal{H}_{0}\right) \\
& = \int_{\varepsilon}^{\infty} p_{0}(x)dx.
\end{aligned}
\end{equation}
For a large $N$, $p_{0}(x)$ can be approximated by a Guassian distribution with mean $\mu_{0}=\sigma_{w}^2$ and variance $\sigma_{0}^2=\frac{1}{N}\sigma_{w}^4$ [12]. Then, the false alarm probability can be expressed as
\begin{equation}
P_f(\varepsilon) = \int_{\varepsilon}^{\infty}\frac{\sqrt{N}}{\sigma_{w}^{2}\sqrt{2\pi}}\mbox{exp}\left(-\frac{N(x-\sigma_{w}^2)^2}{2\sigma_{w}^4}\right)dx.
\end{equation}
Substituting $v=\frac{\sqrt{N}(x-\sigma_{w}^2)}{\sigma_{w}^2}$ into (6),
\begin{equation}
\begin{aligned}
P_f(\varepsilon) & =\frac{1}{\sqrt{2\pi}}\int_{\left(\frac{\varepsilon}{\sigma_{w}^2}-1\right)\sqrt{N}}^{\infty}\mbox{exp}\left(-\frac{v^2}{2}\right)dv\\
& = Q\left(\left(\frac{\varepsilon}{\sigma_{w}^2}-1\right)\sqrt{\tau_{s}f_{s}}\right),
\end{aligned}
\end{equation}
where $Q(x)=\frac{1}{\sqrt{2\pi}}\int_x^{\infty}\mbox{exp}\left(-\frac{u^2}{2}\right)du$ is the complementary distribution function of the standard Gaussian distribution.

Under the hypothesis $\mathcal{H}_{1}$, the detection probability $P_{d}$ defines the probability correctly detecting the presence of primary signal. For a large $N$, $p_{1}(x)$, the PDF of $T(y)$, can be approximated by a Guassian distribution with mean $\mu_{1}=(\gamma_p+1)\sigma_{w}^2$ and variance $\sigma_{1}^2=\frac{1}{N}(\gamma_p+1)^2\sigma_{w}^4$, where $\gamma_{p}={\sigma_{s}^2}/{\sigma_{w}^2}$ is the primary signal-to-noise ratio (SNR) measured at the EH SU [12]. Thus, with the detection threshold $\varepsilon$, 
the detection probability is derived as
\begin{equation}
\begin{aligned}
P_d(\varepsilon) & \triangleq \mbox{Pr}\left(\theta_m=1|\mathcal{H}_{1}\right) \\
& = \mbox{Pr}\left(T(y_{m})>\varepsilon|\mathcal{H}_{1}\right) \\
& = \int_{\varepsilon}^{\infty} p_{1}(x)dx \\
& = Q\left(\left(\frac{\varepsilon}{(\gamma_p+1)\sigma_{w}^2}-1\right)\sqrt{\tau_{s}f_{s}}\right).
\end{aligned}
\end{equation}

Based on the sensing result $\theta_m$, the EH SU can opportunistically access the spectrum. With a high false alarm probability, the idle spectrum is often recognized to be occupied by the PUs. As a result, the improvement of spectrum efficiency is limited, because the EH SU is restricted to use the idle spectrum. However, the EH SU frequently tends to access the occupied spectrum with a low detection probability $P_d$. For highly accurate sensing, the detection probability should be as high as possible, while the false alarm probability should be as low as possible. As mentioned above, the spectrum is considered to be idle at time slot $m$, when the sensing result is $\theta_m=0$. The spectrum access probability for an EH SU thus can be defined as
\begin{equation}
\delta \triangleq \mbox{Pr}(\theta_{m}=0)
\end{equation}
Since there are only two measurable spectrum states, such as idle ($\mathcal{H}_0$) and occupied ($\mathcal{H}_1$), the access probability can be derived according to the {\it{law of total probability}} as follows.
\begin{equation}
\begin{aligned}
\delta & =\sum_{c_{m}\in\{\mathcal{H}_{0}, \mathcal{H}_{1}\}} \mbox{Pr}(\theta_{m}=0|c_{m})\mbox{Pr}(c_{m})\\
&=[1-P_{f}(\varepsilon)]\pi_{i}(q_{i}, q_{o})+[1-P_{d}(\varepsilon)]\pi_{o}(q_{i}, q_{o}).
\end{aligned}
\end{equation}

We assume that the EH SU always has data to transmit. Therefore, the analysis can be applied in the BANs where the biometric information is required to be transmitted in real time.

\section{Analysis of Packet Loss Probability}
We develop a Markovian battery model to capture the energy states of an EH SU equips with an independent EH board and a capacity-limited battery. As illustrated in Fig. 2, the battery is divided into $L$ energy levels with unit of energy $e_u$. We focus on how the energy state changes from one to another depending on the spectrum access and EH. In our work, we assume short sensing duration so that $\tau_s\ll T$. Compared with data transmission, spectrum sensing entails negligible amount of energy. Therefore, energy consumed for spectrum sensing is not considered and the unit of energy $e_u$ denotes sufficient energy for data transmission. In addition, $\chi$ is assumed to be 1 when energy arrival rate is $e_h=\chi e_u$. Suppose that the EH SU always has data to transmit, we analyze transition probabilities at the energy state $l~(0\leq l<L)$ first. Then the energy outage propability that energy is completely ran out ($l=0$) is drawn from the proposed model. Based on sensing inaccuracy and energy outage, we define the packet loss probability.

\subsubsection{$l=0$} If the EH SU runs out of energy, i.e., if the energy state is $l=0$, data transmission is restricted. The energy state will be increased to $l=1$ with probability $p_{0,1}=e_{on}$ as long as the EH board detects the RF signal; otherwise, the energy state will be stuck in $l=0$ with probability $p_{0,0}=1-e_{on}$. During the energy outage period, all of the packets will be doubtlessly lost.
\subsubsection{$0<l<L-1$} Note that from an arbitrary energy state $l$, we can only reach a single unit of energy level, as illustrated in Fig. 2. For example, energy state $l$ can only transit to $l-1$ (lower energy) or $l+1$ (higher energy), or remains in the same state $l$. The energy level will decrease with probability $p_{l, l-1}=\delta (1-e_{on})$ as long as the EH SU successfully detects the idle spectrum when no energy is harvested. However, if the EH board harvests energy while the spectrum is recognized as occupied, the energy level will increase with probability $p_{l, l+1}=(1-\delta) e_{on}$. Otherwise, the energy level $l$ will be maintained with probability $p_{l, l}=(1-\delta)(1-e_{on})+\delta e_{on}$.
\subsubsection{$l=L-1$} Despite much energy being harvested, the energy level $L-1$ will remain unchanged because the battery size is limited. However, the energy level will decrease with probability $p_{l, l-1}=\delta (1-e_{on})$ if no energy is harvested but idle spectrum is detected. Finally, the battery will be fully charged with probability $p_{l, l}=1- \delta(1-e_{on})$.

\begin{figure*}[!t]
\centering
\includegraphics[width=6.5in]{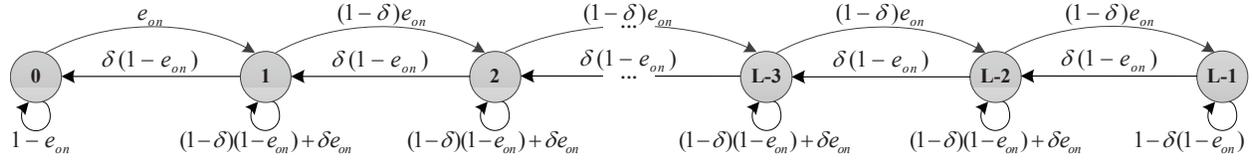}
\caption{A Markovian battery model for the EH SUs}
\label{Fig_sim}
\end{figure*}

Fig. 2 shows the Markovian battery model with the transition probabilities analyzed above. As Fig. 2 illustrates, since all energy states communicate with each other, the proposed finite-state Markovian model is {\it{irreducible}} and therefore {\it{positive recurrent}}. In addtion, starting in arbitrary energy state $l$, we can enter state $l$ at times 1, 2, 4, ..., so that the proposed model is {\it{aperiodic}} because the minimum period is 1. Then, we can derive the steady state probability $\pi_l$ at the energy state $l$ by analyzing the equilibrium equation $\boldsymbol{\pi}=\boldsymbol{\pi}\mathbf{P}$, where $\boldsymbol\pi$ is the vector consisting of steady state probabilities and $\mathbf{P}$ denotes transition matrix and each of them is given as
\begin{equation}
\boldsymbol{\pi}=\left[\begin{array}{cccc}
\pi_0 & \pi_1 & \cdots & \pi_{L-1} \\
\end{array}\right],
\end{equation}

\begin{equation}
\mathbf{P}=\left[\begin{array}{ccccc}
p_{0,0} & p_{0,1} & \cdots & p_{0,L-1} \\
p_{1,0} & p_{1,1} & \cdots & p_{1,L-1} \\
\vdots & \vdots & \cdots & \vdots \\
p_{L-1,0} & p_{L-1,1} & ... & p_{L-1,L-1}\\
\end{array}\right].
\end{equation}

Therefore, at the energy state $l=0$, the equilibrium equation can be written as
\begin{equation}
\pi_{0}e_{on}=\pi_{1}\delta(1-e_{on}).
\end{equation}
\begin{equation}
\pi_{1}=\frac{e_{on}\pi_0}{\delta(1-e_{on})}.
\end{equation}

Furthermore, at the state $l=1$, the analysis can be expressed as
\begin{equation}
\pi_{1}[\delta(1-e_{on})+(1-\delta)e_{on}]=\pi_{0}e_{on}+\pi_{2}\delta(1-e_{on}).
\end{equation}
By using (13), (15) can be simplified as 
\begin{equation}
\begin{aligned}
\pi_{2} & =\frac{(1-\delta)e_{on}}{\delta(1-e_{on})}\pi_{1}\\
& =\left[\frac{(1-\delta)e_{on}}{\delta(1-e_{on})}\right]^2\cdot\frac{\pi_{0}}{(1-\delta)}.
\end{aligned}
\end{equation}

For an arbitrary state $l$ satisfying $1<l<L-1$, we know that the equilibrium equation is given by
\begin{equation}
\pi_{l}[\delta(1-e_{on})+(1-\delta)e_{on}]=\pi_{l-1}(1-\delta)e_{on}+\pi_{l+1}\delta(1-e_{on}).
\end{equation}
When energy state is $l=2$, (17) will become
\begin{equation}
\begin{aligned}
\pi_{3} & =\frac{(1-\delta)e_{on}}{\delta(1-e_{on})}\pi_{2}\\
& =\left[\frac{(1-\delta)e_{on}}{\delta(1-e_{on})}\right]^3\cdot\frac{\pi_{0}}{(1-\delta)}.
\end{aligned}
\end{equation}
To simpify the equation, we define $\alpha$ as
\begin{equation}
\alpha=\frac{(1-\delta)e_{on}}{\delta(1-e_{on})}.
\end{equation}
It is easily known that the similar pattern will be repeated. Thus the analysis for state $L-2$ can be expressed as
\begin{equation}
\pi_{L-1}=\alpha\pi_{L-2}=\alpha^{L-1}\cdot\frac{\pi_{0}}{(1-\delta)}.
\end{equation}
Therefore, we can express $\pi_{l}$ for $1\leq l \leq L-1$ as
\begin{equation}
\pi_{l}=\alpha\pi_{l-1}=\alpha^{l}\cdot\frac{\pi_{0}}{(1-\delta)}.
\end{equation}

Note that the sum of all probabilities is equal to 1, hence $\sum_{l=0}^{L-1} \pi_{l}=1$. Then the equation can be rewrittern as follows.
\begin{equation}
\pi_0+\sum_{l=1}^{L-1}\frac{\alpha^l\pi_0}{\left(1-\delta\right)}=1
\end{equation}
\begin{equation}
\pi_{0}\left[1+\sum_{l=1}^{L-1}\frac{\alpha^l}{(1-\delta)}\right]=1.
\end{equation}
\begin{equation}
\pi_{0}\left[1+\frac{1}{(1-\delta)}\sum_{l=0}^{L-1}\alpha^l-\frac{1}{(1-\delta)}\right]=1.
\end{equation}
Under the constraint that $\alpha<1$,
\begin{equation}
\sum_{l=0}^{L-1}\alpha^l=\frac{1-\alpha^l}{1-\alpha}.
\end{equation}

Substituting (25) into (24), the steady state probability of energy outage depending on the constraint of $\alpha$ is derived as
\begin{equation}
\pi_{0}=\left\{ \begin{array}{rcl}
\cfrac{(1-\delta)(1-\alpha)}{1-\alpha^L-\delta(1-\alpha)} & \mbox{for} & \alpha<1\\
\cfrac{1}{1+\cfrac{1}{(1-\delta)}\sum\limits_{l=1}^{L-1}\alpha^l} & \mbox{for} & \alpha>1.
\end{array}\right.
\end{equation}

If sufficient energy and `real' idle spectrum are provided to EH SUs, secondary data will be reliably transmitted without collisions; otherwise, the packet is most likely to be lost. Then we define the probability of packet loss as
\begin{equation}
\begin{aligned}
P_{L} & \triangleq 1 - \left(1-\pi_{0}\right)\mbox{Pr}\left(\theta_{m}=0, \mathcal{H}_{0}\right) \\
& = 1-\left(1-\pi_{0}\right)\left(1-P_{f}\left(\varepsilon\right)\right)\pi_{i}\left(q_{i}, q_{o}\right).
\end{aligned}
\end{equation}

Wireless channel varies in real-time and therefore channel estimation increases computational complexity. Moreover, estimation errors are inevitable, so as to affect the accuracy of results. Fortunately, faded signals can be normalized by improving the diversity performance of receiver. Therefore, we assume that a single unit of energy guarantees reliable secondary transmission as long as no collision occurs but it can be further extended to a specific wireless channel model.

\section{Simulation and Results}
\begin{table}[!t]
\renewcommand{\arraystretch}{1.2}
\centering
\caption{Simulation Parameters}
\label{my-label}
\begin{tabular}{c|c|c}
\hline
\multicolumn{2}{c|}{Parameter}     & Value      \\
\hhline{===}
\multicolumn{2}{c|}{Battery size, $L$} & 100 \\
\hline
\multicolumn{2}{c|}{Slot duration, $T$} & 100 ms \\
\hline
\multicolumn{2}{c|}{Sensing duration, $\tau_s$} & 2 ms \\
\hline
\multicolumn{2}{c|}{Sampling frequency, $f_s$} & 1 MHz \\
\hline
\multicolumn{2}{c|}{Number of PU} & 1, 10, 20 \\
\hline
\multicolumn{2}{c|}{Number of SU} & 1 \\
\hline
\multirow{6}{*}{Case I} & {Transition probabilities for primary network, ($q_o, q_i$)} & (0.7, 0.5) \\ \cline{2-3}
& \multirow{3}{*}{Transition probabilities for energy harvesting, ($p_{on}, p_{off}$)} & (0.7, 0.5) \\ & & (0.5, 0.5) \\  & & (0.3, 0.5) \\  
\cline{2-3}
& Target false alarm probability, $\overline{P}_f$ & 0.01 \\
\hline
\multirow{6}{*}{Case II} & \multirow{3}{*}{Transition probabilities for primary network, ($q_o, q_i$)} & (0.7, 0.5) \\ & & (0.5, 0.5) \\ & & (0.3, 0.5) \\ \cline{2-3}
& Transition probabilities for energy harvesting, ($p_{on}, p_{off}$) & (0.5, 0.7) \\
\cline{2-3}
& Target primary SNR, $\overline{\gamma}_p$ & -15 dB\\
\hline
\end{tabular}
\end{table}

Monte-Carlo simulation is performed under different processes of energy arrival and primary network, and numerical results are used to evaluate the proposed analysis. We mainly focus on how the packet loss probability changes with the variations of primary SNR $\gamma_p$ (Case I) and normalized detection threshold $\epsilon/\sigma_{w}^2$ (Case II). Table 1 provides the parameters utilized throughout the simulation.

Figs. 3-4 illustrate the packet loss probability as a function of the primary SNR $\gamma_p$ for \{$q_o=0.7, q_i=0.5$\} and \{$p_{on}=0.5, p_{off}=0.7$\}, respectively. We set the target false alarm probability to $\overline{P}_f=0.01$. Note that higher primary SNR improves reliability of spectrum sensing ($P_d:0\rightarrow1$). Therefore, miss detection rarely happens and more energy can be accumulated at the EH SU with higher primary SNR. As a result, the EH SU rarely runs out of energy the packet loss probability decreases accordingly as the primary SNR increases.

\begin{figure} [!t]
\centering
\includegraphics[width=4in]{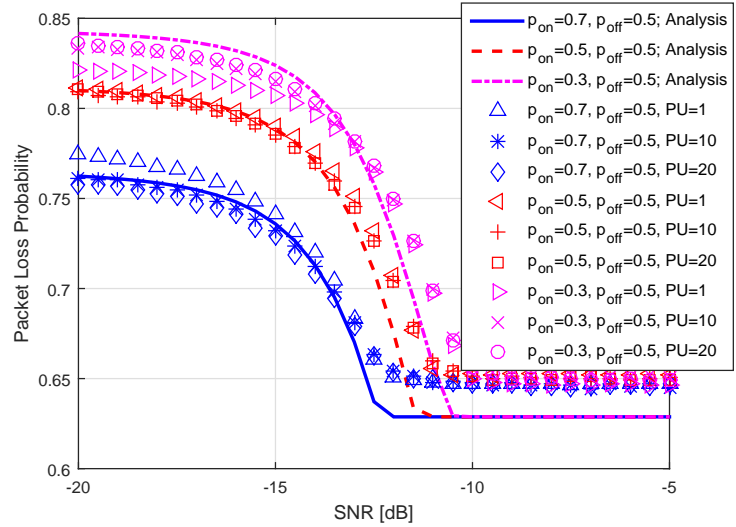}
\caption{Packet loss probability vs. primary SNR $\gamma_p$ for $P_f=0.01$, $q_o=0.7$ and $q_i=0.5$.}
\label{Fig_sim}
\end{figure}
\begin{figure} [!t]
\centering
\includegraphics[width=4in]{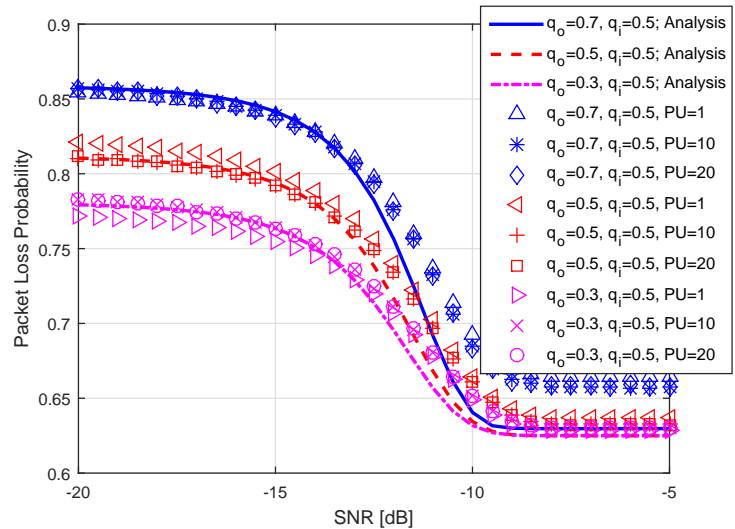}
\caption{Packet loss probability vs. primary SNR $\gamma_p$ for $P_f=0.01$, $p_{on}=0.5$ and $p_{off}=0.7$.}
\label{Fig_sim}
\end{figure}

\begin{figure} [!t]
\centering
\includegraphics[width=4in]{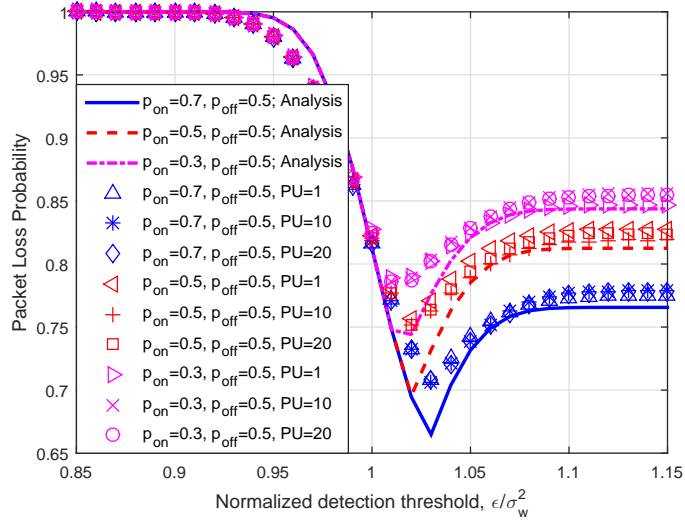}
\caption{Packet loss probability vs. normalized detection threshold $\varepsilon/\sigma_{w}^2$ for $\gamma_p=-15~{\rm{dB}}$, $q_o=0.7$ and $q_i=0.5$.}
\label{Fig_sim}
\end{figure}
\begin{figure} [!t]
\centering
\includegraphics[width=4in]{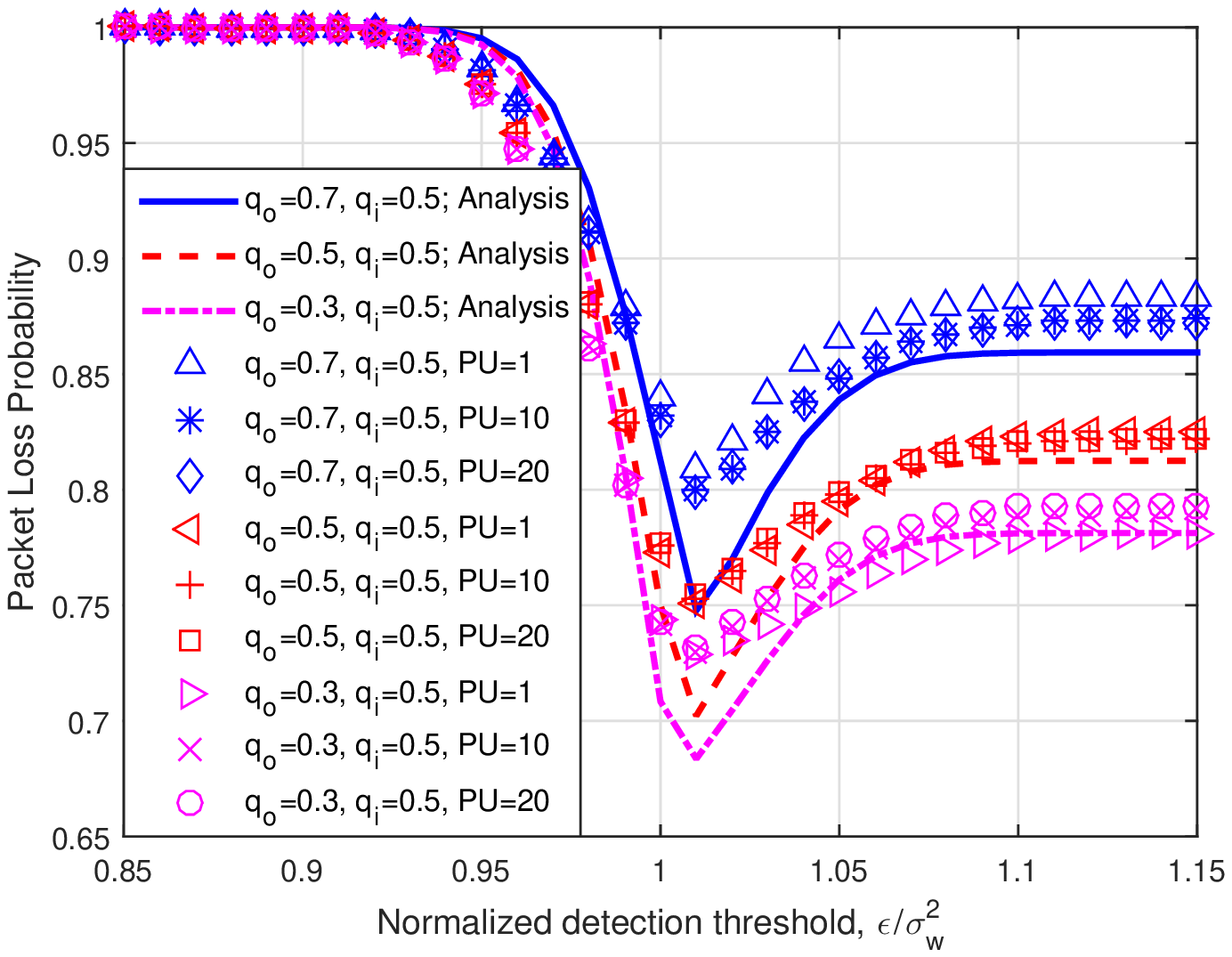}
\caption{Packet loss probability vs. normalized detection threshold $\varepsilon/\sigma_{w}^2$ for $\gamma_p=-15~{\rm{dB}}$, $p_{on}=0.5$ and $p_{off}=0.7$.}
\label{Fig_sim}
\end{figure}

Figs. 5-6 present the packet loss probability corresponding to the normalized detection threshold $\varepsilon/\sigma_{w}^2$ that indicates the relative detection threshold compared to the noise power. Similarly, \{$q_o=0.7, q_i=0.5$\} and \{$p_{on}=0.5, p_{off}=0.7$\} are considered, and the primary SNR is $\gamma_p=-15$ dB. As the normalized detection threshold increases, false alarm is improved ($P_f:1\rightarrow0$), while miss detection happens more frequently ($P_d:1\rightarrow0$). At the beginning, the number of secondary data packets lost is decreasing, since the access probability is increasing. However, more frequent spectrum accessing results in energy outage and collisions. Consequently, increasing the number of secondary data packets does not return a gain when the normalized detection threshold keeps increasing.

Figs. 3-6 also show that the minimum and the maximum probabilities are achieved with frequent energy harvesting and occasional spectrum accessing, respectively. In addition, no matter how many PUs are involved to share resources with the EH SU, secondary packets are lost due to the energy outage and the sensing inaccuracy since the EH SU cannot sense all spectrums simultaneously and selects one of spectrums to sense. Finally, we explain the reason why the packet loss probability becomes saturated, as follows. The access probability converges to a constant as $P_d\rightarrow1$ when $P_f$ is fixed or both $P_f$ and $P_d$ simultaneously approach 0, so that the energy outage probability becomes saturated.

\section{Conclusions}
To derive the probability of packet loss in the EH SUs, we have proposed a Markovian battery model. Form the battery model, the energy outage probability was analyzed first. Then the packet loss probability due to sensing inaccuracy and energy outage was derived. Through the Monte-Carlo simulation, the analysis was shown to be acceptable to predict the packet loss probability and therefore can be applied to the cross layer optimization.

\end{document}